\documentclass[reqno, oneside]{amsart}
\usepackage{style}

\begin{document}
\title[Muckenhoupt Hamiltonians]{Sampling measures, Muckenhoupt Hamiltonians, and triangular factorization}
\author{R.~V.~Bessonov}

\address{St.Petersburg State University ({\normalfont \hbox{7/9}, Universitetskaya nab., 
St.Petersburg, 199034 Russia}) and St.Petersburg Department of Steklov Mathematical Institute of Russian Academy of Science ({\normalfont 27, Fon\-tan\-ka, St.Petersburg, 191023 Russia})}
\email{bessonov@pdmi.ras.ru}

\thanks{The work is supported by RFBR grant mol\_a\_dk 16-31-60053 and by ``Native towns'', a social investment program of PJSC "Gazprom Neft"}

\begin{abstract}
Let $\mu$ be an even measure on the real line $\R$ such that 
$$c_1 \int_\R|f|^2\,dx \le \int_\R|f|^2\,d\mu \le c_2\int_\R|f|^2\,dx$$ 
for all functions $f$ in the Paley-Wiener space $\mathrm{PW}_{a}$. We prove that $\mu$ is the spectral measure for the unique Hamiltonian $\H=\left(\begin{smallmatrix}w&0\\0&\frac{1}{w}\end{smallmatrix}\right)$ on $[0,a]$ generated by a weight $w$ from the Muckenhoupt class $A_2[0,a]$. As a consequence of this result, we construct Krein's  orthogonal entire functions with respect to $\mu$ and prove that every positive, bounded, invertible Wiener-Hopf operator on $[0,a]$ with real symbol admits triangular factorization.
\end{abstract}

\maketitle

\section{Introduction}\label{s0}
The classical Paley-Wiener space $\pw_a$ consists of entire functions of exponential type at most $a$ square summable on the real line, $\R$. A measure $\mu$ on $\R$ is called a sampling measure for the space $\pw_a$ if there exist positive constants $c_1$, $c_2$ such that 
\begin{equation}\label{eq1}
c_1 \int_\R|f|^2\,dx \le \int_\R|f|^2\,d\mu \le c_2\int_\R|f|^2\,dx, \qquad f \in \pw_a.
\end{equation}
Let $\H$ be a regular Hamiltonian on $[0,a]$, that is, $\H$ is a mapping from $[0,a]$ to the set of $2\times 2$ non-negative matrices with real entries such that $\trace\H$ is a positive non-vanishing function in $L^1[0, a]$. Denote by $\Theta_\H = \Theta_\H(r,z)$ solution of the following Cauchy problem:
\begin{equation}\label{eq2}
JX'(r) = z\H(r)X(r), \quad X\colon [0, a] \to \C^2, \quad X(0) = \oz, \quad z \in \C.  
\end{equation}
It is known from a general theory of canonical Hamiltonian systems that for every measure $\mu$ satisfying \eqref{eq1} there exists a regular Hamiltonian $\H$ with $\int_{0}^{a}\sqrt{\det\H} = a$ such that $\mu$ is a spectral measure for problem \eqref{eq2}. The latter means that the Weyl-Titch\-marsh transform 
\begin{equation}\label{eq3}
\W_{\H,a}\colon X \mapsto \frac{1}{\sqrt{\pi}}\int_{0}^{a} \bigl\langle\H(r) X(r),\Theta_\H(r, \bar z)\bigr\rangle_{\C^2}\,dr, \qquad z \in \C,
\end{equation}
generated by solution $\Theta_\H$ of Cauchy problem \eqref{eq2} maps isometrically the space 
\begin{gather}
L^2(\H, a) = \Bigl\{X\!:\![0, a] \to \C^2 \colon \|X\|_{L^2(\H,a)}^{2} \! = \! \int_{0}^{a}\!\bigl\langle\H(r) X(r),X(r)\bigr\rangle_{\C^2}dr < \infty \Bigr\}\Big/\!{\mathcal K}(\H), \notag\\
{\mathcal K}(\H) = \Bigl\{X\colon \;\H(t) X(t) = 0 \mbox{ for almost all } t\in [0,r]\Bigl\} \notag
\end{gather}
into the space $L^2(\mu)$. A general problem in the inverse spectral theory is to translate properties of a spectral measure $\mu$ into properties of the Hamiltonian $\H$ it generates. 

\medskip

Two essentially different cases of the above problem attracted much attention. If $\mu$ is a ``small perturbation'' of the Lebesgue measure on $\R$ (in the sense that the Fourier transform of $\mu$ restricted to the interval $[-a,a]$ differs from the point mass measure $\delta_0$ concentrated at $0$ by a function in $L^1[-a,a]$), the I.~M.~Gelfand--B.~M.~Levitan approach \cite{GL55}, \cite{March11} gives a quite precise information on relation between $\mu$ and $\H$. On the other hand, if $\mu$ is arbitrary measure on $\R$ such that $\int_\R\frac{d\mu(t)}{1+t^2} < \infty$, the theory of M.~G.~Krein \cite{Krein74} (for even measures $\mu$) and  L.~de~Bran\-ges~\cite{dbbook} (for all~$\mu$) implies the {\it existence} of a unique Hamiltonian $\H \in L^1_{{\rm loc}}[0, \infty)$ such that $\mu$ is the spectral measure for $\H$. However, it is not known how translate even simple properties of a Hamiltonian $\H$ (e.g., membership in $L^p$ class for some $p>1$) to the properties of its spectral measure $\mu$ and vice versa.  In this paper we consider a ``median'' situation (spectral measures with sampling property \eqref{eq1} for the Paley-Wiener space $\pw_a$) and use both Gelfand-Levitan and Krein-de Branges theories.

\medskip

A measure $\mu$ on $\R$ is called even if $\mu(S) = \mu(-S)$ for every Borel set $S\subset \R$. A function $w>0$ belongs to the Muckenhoupt class $A_2[0,a]$ if the supremum of products $\bigl(\frac{1}{|I|}\int_{I}w\bigr)\!\cdot\!\bigl(\frac{1}{|I|}\int_{I}\frac{1}{w}\bigr)$ over all intervals $I \subset [0,a]$ is finite. Here is the main result of the paper. 
\begin{Thm}\label{t1} 
Let $\mu$ be an even sampling measure for $\pw_a$. Then $\mu$ is the spectral measure for problem \eqref{eq2} corresponding to the unique Hamil\-tonian $\H = \left(\begin{smallmatrix}w&0\\0&\frac{1}{w}\end{smallmatrix}\right)$ generated by a weight $w \in A_2[0,a]$.
\end{Thm}
The Hamiltonian $\H$ in Theorem \ref{t1} could be recovered from the spectral measure~$\mu$ by means of the following simple formula:
$$
w(r) = \pi\frac{\partial}{\partial r} \Bigl\|T_{\mu,r}^{-1}\sinc_{r}\Bigr\|^{2}_{L^2(\mu)}, \qquad
\sinc_{r} = \frac{\sin rx}{\pi x}, \qquad r \in [0,a],
$$
where $T_{\mu,r}$ is the truncated Toeplitz operator on $\pw_r$ with symbol $\mu$ defined by  
\begin{equation}\label{eq4}
(T_{\mu,r} f)(z) = \int_{\R}f(x) \frac{\sin r(x - z)}{\pi(x - z)}\,d\mu(x), \qquad z \in \C.
\end{equation}
A nontrivial fact is that the continuous increasing function $r \mapsto \|T_{\mu,r}^{-1}\sinc_{r}\|^{2}_{L^2(\mu)}$ is {\it absolutely} continuous and its derivative $w/\pi$ does not vanish on a set of positive Lebesgue measure. In the proof of Theorem \ref{t1} we first obtain an estimate for the ``$A_2$-norm'' of $w$ in terms of $c_1$, $c_2$ assuming above properties of $w$; then use an approximation argument based on a description of positive truncated Toeplitz operators on $\pw_r$ and $L^p$-summabilty of weights $w \in A_2[0,a]$ for some $p>1$. 

\medskip

Section 5 in \cite{BR2015} contains an example of a diagonal Hamiltonian $\H$ on $[0, 1]$ such that both $\H$, $\H^{-1}$ are uniformly bounded on $[0,1]$, but the spectral measures of the corresponding problem \eqref{eq2} fail to have sampling property. This shows that $A_2[0,a]$ class does not describe canonical Hamiltonian systems generated by sampling measures for $\pw_a$.

\medskip

Theorem \ref{t1} yields two results of independent interest. 

\medskip

Given a measure $\mu$ satisfying \eqref{eq1} and a number $r \in [0,2a]$, denote by $(\pw_{[0,r]},\mu)$ the Paley-Wiener space of functions from $L^2(\R)$ with Fourier spectrum in $[0,r]$ equipped with the inner product taken from $L^2(\mu)$.

\begin{Thm}\label{t2}
Let $\mu$ be an even sampling measure for the space $\pw_a$. Then there exists a family of entire functions $\{P_{t}\}_{t \in [0, 2a]}$ such that $\F_\mu: f \mapsto \frac{1}{\sqrt{2\pi}}\int_{0}^{r}f(t)P_t(z)\,dt$ is the unitary operator from $L^2[0,r]$ to $(\pw_{[0,r]},\mu)$ for every $r \in [0,2a]$.
\end{Thm}
In the case where $\mu$ is a ``small perturbation'' of the Lebesgue measure (see discussion above), the functions $P_r$ in Corollary \ref{t2} coincide with orthogonal entire functions constructed by M. G. Krein in \cite{Kr81}. S. A. Denisov provides an extensive treatment of the subject, collecting many old and new results in paper \cite{Den06}. 

\medskip

The second application of Theorem \ref{t1} concerns the classical factorization problem for positive invertible operators. Let $H$ be a separable Hilbert space and let $B(H)$ be the algebra of all bounded operators on $H$. Consider a complete chain $\N$ of subspaces in $H$ and denote by $\A_\N = \{A \in B(H): A E\subset E, \; E \in \N\}$ the nest algebra of upper-triangular operators with respect to $\N$. 
In sixties, I.~C.~Gohberg and M.~G.~Krein proved (see Theorem 6.2 in Chapter 4 of \cite{GK70}) that every positive invertible operator $T$ on $H$ of the form $T=I-K$ with $K$ in Macaev ideal $S_{\omega}$ admits the triangular factorization $T = A^{\ast}A$, where $A= I-K_A$ is an invertible operator on $H$ such that $K_{A} \in S_\omega \cap \A_\N$. Famous theorem by D. R. Larson \cite{Larson85} says that every positive invertible operator $T$ admits triangular factorization $T = A^*A$ with $A, A^{-1} \in \A_\N$ if and only if the chain $\N$ is countable. Moreover, given $0<\eps<1$, the non-factorable operator $T$ can be chosen so that $K = I - T$ is a compact operator with  
$\|K\| < \eps$.  

We consider the problem of triangular factorization for Wiener-Hopf convolution operators. Let $\psi \in \Sch'$ be a tempered distribution on $\R$ and let $0<a\le \infty$. The Wiener-Hopf operator $W_{\psi}$ on $L^2[0,a)$ with symbol $\psi$ is densely defined by
$$
(W_{\psi}f)(y) = \bigl\langle\psi, s_y f\bigr\rangle_{\Sch'}, \qquad  y \in [0, a), \qquad s_y f: x \mapsto f(x - y),
$$ 
on smooth functions $f$ with compact support in $(0,a)$. In the case where $\psi\in L^1(\R)$ we have more familiar definition, $W_{\psi}\colon f \mapsto \int_{0}^{a}\psi(x-y)f(x)\,dx$. As following result shows, Wiener-Hopf operators with real symbols are always factorable.
\begin{Thm}\label{t3}
Let $0<a \le \infty$. Every positive, bounded, and  invertible Wiener-Hopf operator $W_\psi$ on $L^2[0,a)$ with real symbol $\psi \in \Sch'$ admits triangular factorization: $W_\psi = A^\ast A$, where $A$ is a bounded invertible operator such that $A L^2[0,r] = L^2[0,r]$ for every $r \in [0,a)$.
\end{Thm}
Wiener-Hopf operators $W_\psi$ in Theorem \ref{t3} admit triangular factorizations in the reverse order $W_\psi =  AA^\ast$ as well. Relation of absolute continuity of aforementioned function $r \mapsto \|T_{\mu,r}^{-1}\sinc_{r}\|_{L^{2}(\mu)}^{2}$ to triangular factorization problems has been previously found in different terms by L. A. Sakhnovich, see Theorem 4.2 in \cite{Sakhnovich07}. On the other hand, Theorem \ref{t3} contradicts Theorem 4.1 from another work \cite{Sakhnovich12} by the same author. See discussion in Section \ref{s5}.

\medskip

\noindent {\bf Acknowledgement.} The author is grateful to many colleagues who took a part in discussions related to the subject of the paper, especially to Roman Romanov, Mikhail Sodin, Pavel Zatitsky, and Dmitriy Zaporozhets. \\


\section{Integration over simplex and the Muckenhoupt class $A_2$}\label{s1}
\noindent Let $w$ be a positive function on an interval $[0,a]$. We associate to $w$ the quantity 
$$
\|w\|_{A_2[0,a]} = \sup_{I\subset [0,a]}\left(\frac{1}{|I|}\int_{I}w(x)\,dx\right)\!\cdot\!\left(\frac{1}{|I|}\int_{I}\frac{1}{w(x)}\,dx\right),
$$
where $I$ runs over all subintervals of $[0,a]$. Note that $\|\cdot\|_{A_2[0,a]}$ is not a norm in the standard sense, but we will use this convenient notation. The Muckenhoupt class $A_2[0,a]$ consists of functions $w>0$ such that $\|w\|_{A_2[0,a]} < \infty$. In this section we present a special integral condition for a weight $w$ to belong to the $A_2[0,a]$ class.

\medskip

Let $\phi$ be a real-valued function on the interval $[0, a]$. For a real $0<t<a$ and an integer $n \ge 1$ define the mapping 
\begin{equation}\label{eq5}
G_{\phi,n}: x \mapsto \sum_{k=1}^{n} (-1)^{n+k}\phi(x_k), \quad x \in K_{t,n}, 
\end{equation}
on simplex   
$
K_{t,n} = \{x \in \R^n: x = (x_1, \ldots, x_n), \; t\ge x_1 \ge \ldots \ge x_n \ge 0\}. 
$ 
Let $m_n$ denote the usual Lebesgue measure on $\R^n$.

\medskip

Next proposition will be used in the proof of Theorem \ref{t1}.
\begin{Prop}\label{p1}
Let $\phi$ be a function on $[0, a]$ such that $e^{|\phi|} \in L^1[0,a]$. Assume that for every $r \in [0, a]$ and every integer $n \ge 1$ we have
\begin{gather}
\frac{1}{a_n(r)}\int_{0}^{r}e^{(-1)^n \phi(t)}\left(\int_{K_{t,n}}e^{G_{\phi,n}(x)}dm_n(x)\right)^2 dt \le b_2, \label{eq6} \\ 
b_1 \le \frac{1}{r}\int_{0}^{r}e^{\phi(t)}\,dt \le b_2; \label{eq7}
\end{gather}
where $b_1$, $b_2$ are positive constants, and $a_n(r) = r^{2n+1}(2n+1)^{-1}(n!)^{-2}$. 
Then the function $w = e^{\phi}$ belongs to $A_2[0,a]$ and $\|w\|_{A_2[0,a]} \le 2^{28}(b_2 + b_{1}^{-2}b_2)^{14}$. 
\end{Prop}

\medskip

We first prove several preliminary estimates.

\begin{Lem}\label{l1}
Let $\phi$ be a function as in Proposition \ref{p1}. Then for every $r \in [0,a]$ and $b = 2(b_2 + b_1^{-2}b_2)$ we have  	
\begin{equation}\label{eq8}
\frac{1}{r}\int_{0}^{r}|\phi(t)|\,dt \le \log b, \qquad \frac{1}{r}\int_{0}^{r}e^{|\phi(t)|}\,dt \le b.
\end{equation}
Consequently, for every decreasing differentiable function $k\ge 0$ on $[0,r]$ satisfying $\int_{0}^{r}k(t)\,dt = 1$ and 
$k(r) = 0$ we have $\int_{0}^{r}|\phi(t)|k(t)\,dt \le \log b$.
\end{Lem}
\beginpf Clearly, the first estimate in \eqref{eq8} follows from the second one and the Jensen's inequality for convex function $e^x$. Taking $n = 1$ in \eqref{eq6}, we obtain
$$
\frac{3}{r^3}\int_{0}^{r}e^{-\phi(t)}\left(\int_{0}^{t}e^{\phi(t_1)}\,dt_1\right)^2\,dt \le b_2. 
$$
From \eqref{eq7} we know that $\frac{1}{t}\int_{0}^{t}e^{\phi(t_1)}\,dt_1 \ge b_1$  for all $t \in [0,r]$. It follows that
$$
b_1^{-2}b_2 \ge \frac{3}{r^3}\int_{0}^{r}e^{-\phi(t)}t^2\,dt \ge \frac{1}{r}\int_{r/2}^{r}e^{-\phi(t)}\,dt.
$$ 
Using the other side estimate $\frac{1}{r}\int_{0}^{r}e^{\phi(t)}\,dt \le b_2$ and inequality $e^{|x|} \le e^x + e^{-x}$, we see that
$$
\frac{2}{r}\int_{r/2}^{r}e^{|\phi(t)|}\,dt \le b_2 + b_1^{-2}b_2
$$
for all $r \in [0,a]$. Then \eqref{eq8} follows from
$$
\frac{1}{r}\int_{0}^{r}e^{|\phi(t)|}\,dt = \frac{1}{r}\left(\sum_{k=0}^{\infty}|I_{r,k}|\cdot\frac{1}{|I_{r,k}|}\int_{I_{r,k}}e^{|\phi(t)|}\,dt \right) \le b, 
$$
where $I_{r,k} = [2^{-k-1}r, 2^{-k}r]$. Now if $k$ is a function on $[0,r] \subset[0,a]$ as in the statement, we have
\begin{align*}
\int_{0}^{r}|\phi(t)|k(t)\,dt 
&= -\int_{0}^{r}|\phi(t)|\int_{0}^{r}\chi_{[t,r]}(s)k'(s)\,ds\,dt \\
&= -\int_{0}^{r}k'(s)\int_{0}^{r}\chi_{[0,s]}(t)|\phi(t)|\,dt\,ds \\
&\le - \log b \int_{0}^{r}k'(s)s\,ds = \log b.   
\end{align*}
This completes the proof. \qed

\medskip

For $n\ge 1$ introduce the intervals $I_{t,n} = [\delta_n t, t]$, where $\delta_n = 1 - \frac{1}{n+1}$ if $n$ is odd, and $\delta_n = 1 - \frac{1}{n}$ if $n$ is even.
In particular, $I_{t,n} = I_{t,n+1}$ for every odd $n$. Set 
$$
[\phi]_{t,n} = 2(-1)^{n+1}\int_{K_{t,n}}G_{\phi,n}(x)\,dm_{t,n}(x),
$$
where $m_{t, n} = \frac{n!}{t^n}\cdot m_n$ is the scalar multiple of the Lebesgue measure $m_n$ on $\R^n$ normalized so that $m_{t,n}(K_{t,n}) = 1$. 

\medskip

\begin{Lem}\label{l2}
For $r \in [0,a]$ and odd $n \ge 1$ we have $|[\phi]_{\delta_nr,n} - [\phi]_{\delta_{n+1}r,n+1}| < 6\log b$, where $b$ is the constant from Lemma \ref{l1}. 
\end{Lem}
\beginpf Arguing by induction, it is easy check that for all $n \ge 1$ and $\tau \in [0,a]$ we have
$$
[\phi]_{\tau,n} = \int_{0}^{\tau}\phi(s)k_{\tau,n}(s)\,ds, \qquad k_{\tau,n}(s) = \frac{2n}{\tau^n}(2s-\tau)^{n-1}.
$$
For odd (correspondingly, even) integers $n$ the kernels $k_{\tau,n}$ are even (correspondingly, odd) functions with respect to the point $\tau/2$. As $n$ tends to infinity, the kernels $k_{\tau,n}$ tend to zero uniformly on every closed interval in $(0,\tau)$. We also have
\begin{equation}\label{eq9}
\int_{0}^{\tau}|k_{\tau,n}(s)|\,ds = 2, \qquad \sup_{s \in [\tau/2, \tau]}|k_{\tau,n}(s) - k_{\tau,n+1}(s)| \le \frac{2}{\tau}.
\end{equation}
Now take an odd integer $n \ge 1$ and note that $\delta_{n} = \delta_{n+1} = \frac{1}{n+1}$. Setting $\tau = \delta_n r$, we obtain
\begin{align*}
\Bigl|[\phi]_{\tau,n} - [\phi]_{\tau,n+1}\Bigr| 
\le &\int_{0}^{\tau/2}|\phi(s)k_{\tau,n}(s)|\,ds  \\
&+\int_{0}^{\tau/2}|\phi(s)k_{\tau,n+1}(s)|\,ds\\
&+\int_{\tau/2}^{\tau}|\phi(s)|\cdot|k_{\tau,n}(s) - k_{\tau,n+1}(s)|\,ds.
\end{align*}
By Lemma \ref{l1} for functions $\phi$, $k = \frac{1}{2}|k_{\tau,n}|$, and $k = \frac{1}{2}|k_{\tau,n+1}|$ on $[0,\frac{\tau}{2}]$, the sum of first two integrals is bounded from above by $4\log b$. To show that the last integral does not exceed $2\log b$, use \eqref{eq8} and the second estimate in \eqref{eq9}. \qed

\medskip

\noindent{\bf Proof of Proposition \ref{p1}.}  Take an odd integer $n \ge 1$. Since the integrand in \eqref{eq6} is positive, we have
\begin{align*}
b_2 
&\ge \frac{1}{a_n(r)}\int_{\delta_n r}^{r}e^{-\phi(t)}\left(\int_{K_{t,n}}e^{G_{\phi,n}(x)}dm_n(x)\right)^2 dt,\\
&\ge \frac{1}{a_n(r)}\left(\int_{\delta_n r}^{r}e^{-\phi(t)}\,dt\right) \cdot \left(\int_{K_{\delta_nr,n}}e^{G_{\phi,n}(x)}dm_n(x)\right)^2.
\end{align*}
By Jensen's inequality, 
$$
\int_{K_{\delta_nr,n}}e^{G_{\phi,n}(x)}dm_n(x) \ge \frac{(\delta_n r)^n}{n!}\exp\left(\frac{[\phi]_{\delta_n\!r,n}}{2}\right).
$$
For all $n \ge 1$ we have  
$$
\frac{1}{a_n(r)}\cdot \left(\frac{(\delta_nr)^{n}}{n!}\right)^{2} = \frac{(2n+1)(n!)^{2}}{r^{2n+1}}\cdot\frac{r^{2n}}{(n!)^2}\delta_{n}^{2n} \ge \frac{n+1}{32r} = \frac{1}{32|I_{n,r}|}.
$$
We now see that 
\begin{equation}\label{eq10}
\frac{1}{|I_{r, n}|}\int_{I_{r,n}}\exp\Bigl(-\phi(t) + [\phi]_{\delta_nr,n}\Bigr)\,dt \le 32b_2.
\end{equation}
Analogously, for the even integer $n+1$ we have
\begin{equation}\notag
\frac{1}{|I_{r,n+1}|}\int_{I_{r,n+1}} \exp\Bigl(\phi(t) - [\phi]_{\delta_{n+1} \!r, n}\Bigr) \,dt \le 32b_2. 
\end{equation}
Recall that $I_{t,n+1} = I_{t,n}$. Applying Lemma \ref{l2}, we obtain 
\begin{equation}\label{eq11}
\frac{1}{|I_{r,n}|}\int_{I_{r,n}} \exp\Bigl(\phi(t) - [\phi]_{\delta_{n} \!r, n}\Bigr) \,dt \le 32b_2e^{6\log b} \le 32b^7,
\end{equation}  
where $b$ is the constant from Lemma \ref{l1}. Using inequality $e^{|x|} \le e^x + e^{-x}$, we get from \eqref{eq10} and \eqref{eq11} the estimate
\begin{equation}\label{eq12}
\frac{1}{|I|}\int_{I} e^{|\phi(t) - c_I|}\,dt \le 64b^{7}
\end{equation} 
for all intervals $I$ of the form $I = [(1-\frac{1}{n+1})r, r]$, where $r \in [0,a]$, and integer $n \ge 1$ is odd. Here $c_I$ is a constant depending on $I$ (in fact, $c_I = [\phi]_{\delta_n r,n}$ works, but from now on the particular choice of $c_I$ plays no role). Formula \eqref{eq8} gives \eqref{eq12} with $c_I = 0$ for intervals of the form $I = [0, t]$. 

\medskip

Next, observe that each interval $J \subset [0,a]$ is contained in an interval $I$ satisfying~\eqref{eq12} and such that $|I| \le 2|J|$. Indeed, let $t$ be the right point of $J$. If $|J| \ge |t|/2$, take $I = [0,t]$. In the case $|J| < |t|/2$ find an odd number $n\ge1$ such that $I_{t,n+2} \subset J \subset I_{t, n}$ and take $I = I_{t, n}$. Fix this interval $I$ and the corresponding constant $c_I$ form \eqref{eq12}. We have
\begin{align*}
\left(\frac{1}{|J|}\int_{J} e^{\phi}\,dt\right)\!\cdot\!\left(\frac{1}{|J|}\int_{J} e^{-\phi}\,dt\right) 
&\le \left(\frac{2}{|I|}\int_{I} e^{\phi}\,dt\right)\!\cdot\!\left(\frac{2}{|I|}\int_{I} e^{-\phi}\,dt\right) \\
&\le \left(\frac{2}{|I|}\int_{I} e^{\phi-c_I}\,dt\right)\!\cdot\!\left(\frac{2}{|I|}\int_{I} e^{-\phi+c_I}\,dt\right) \le (2b)^{14}.
\end{align*}
Since interval $J$ is arbitrary, this shows that function $w = e^{\phi}$ belongs to the Mucken\-houpt class $A_2[0,a]$ and $\|w\|_{A_2[0,a]} \le (2b)^{14} = 2^{28}(b_2 + b_{1}^{-2}b_2)^{14}$.
\qed


\medskip

\section{Proof of Theorem \ref{t1}}\label{s2}
As it was mentioned in the Introduction, we will use an approximation argument in the proof of Theorem \ref{t1}. To have a stable approximation, we need a result describing positive truncated Toeplitz operators on $\pw_a$.

\medskip

\subsection{Preliminaries on truncated Toeplitz operators. } Let $\mu \ge 0$ be a measure on the real line $\R$ such that $\|f\|_{L^2(\mu)}^{2} \le c \|f\|_{L^2(\R)}^{2}$ for all functions $f \in \pw_{[0, a]}$. Define the truncated Toeplitz operator $A_{\mu,a}$ on $\pw_{[0,a]}$ by the sesquilinear form
\begin{equation}\label{eq13}
(A_{\mu,a} f, g)_{L^2(\R)} = \int_{\R}f\bar g\,d\mu, \qquad f,g \in \pw_{[0,a]}.
\end{equation}
In the case where $\mu = \gphi\,dm$ is absolutely continuous with respect to the Lebesgue measure $m$ on $\R$ and has density $\gphi$, the operator $A_{\mu,a}$ coincides with the projection of the standard Toeplitz operator $T_\gphi$ on the Hardy space $H^2$ to the  subspace $\pw_{[0,a]}$. This explains the name ``truncated Toeplitz'' for the operator $A_{\mu,a}$. 

\medskip

It is well-known (see, e.g., Section 6.1 in \cite{Nik02t1}) that the operator 
\begin{equation}\label{eq14}
V: h \mapsto \frac{1}{\sqrt{\pi}}\frac{1}{z+i}h\left(\frac{z-i}{z+i}\right), \quad z \in \C_+,
\end{equation}
maps unitarily the Hardy space $H^2(\D)$ in the open unit disk $\D = \{\xi \in \C: |\xi|<1\}$ onto the Hardy space $H^2$ in the upper half-plane $\C_+ = \{z \in \C: \Im z >0\}$. Moreover, for every $a>0$ we have $VK_{\theta_{a}} = \pw_{[0,a]}$, where $\theta_{a} = \exp\bigl(a\frac{z+1}{z-1}\bigr)$ is the inner function in $\D$ and $K_{\theta_{a}}$ is the orthogonal complement in $H^2(\D)$ to the subspace $\theta_{a} H^2(\D)$. As we will see in a moment, the truncated Toeplitz operators defined by~\eqref{eq13} are unitarily equivalent to truncated Toeplitz operators on the shift-coinvariant subspace $K_{\theta_a}$ of $H^2(\D)$. See D. Sarason's paper \cite{Sar07} for basic properties of truncated Toeplitz operators on general coinvariant subspaces of $H^2(\D)$. 

\medskip 

We also will  deal with the operators $T_{\mu,a}$ on the space $\pw_a$ defined by the same sesquilinear form 
$$
(T_{\mu,a} f, g) = \int_{\R}f\bar g\,d\mu, \qquad f,g \in \pw_a.
$$
It is easy to see that this definition agrees with formula \eqref{eq4}. By construction, we have $T_{\mu, a} = V_{a}^{-1} A_{\mu, 2a} V_{a}$, where $V_a : \pw_a \to \pw_{[0,2a]}$ is the unitary operator taking a function $f$ into $e^{iaz} f$.

\medskip

\begin{Lem}\label{l3}
Let $T$ be a positive bounded operator on $\pw_{[0,a]}$ satisfying relation 
\begin{equation}\label{eq15}
(Tf,f)_{L^2(\R)} = (T\tfrac{z-i}{z+i}f,\tfrac{z-i}{z+i}f)_{L^2(\R)}  
\end{equation}
for all functions $f \in \pw_{[0,a]}$ such that $f(-i) = 0$. Then there exists a positive measure $\mu$ on $\R$ such that $T=A_{\mu,a}$. Similarly if $T$ is a positive bounded operator on $\pw_a$ satisfying \eqref{eq15} for all $f \in \pw_{a}$ such that $f(-i) = 0$, then $T = T_{\mu,a}$ for a positive measure $\mu$ on~$\R$. 
\end{Lem}
\beginpf Let $\theta_a$, $K_{\theta_a}$, and $V: K_{\theta_a} \to \pw_{[0,a]}$ be defined as above. Consider the operator $\tilde T = V^{-1}TV$ on $K_{\theta_{a}}$ unitarily equivalent to the operator $T$ on $\pw_{[0,a]}$. Recall that the inner product in $K_{\theta_a}$ is inherited from the space $L^2(\T)$ on the unit circle $\T = \{\xi\in\C:\;|\xi| =1\}$. Assumption \eqref{eq15} means that 
\begin{equation}\label{eq16}
(\tilde T h, h)_{L^2(\T)} = (\tilde T \xi h, \xi h)_{L^2(\T)}
\end{equation}
for every function $h \in K_{\theta_{a}}$ such that $\xi h \in K_{\theta_{a}}$. Indeed,  $(V\xi h)(z) = \tfrac{z-i}{z+i}(Vh)(z)$ and hence $V(\xi h) \in \pw_{[0,a]}$ if and only if $(Vh)(-i) = 0$. Theorem 8.1 in \cite{Sar07} says that a bounded operator $\Tilde T$ on $K_{\theta_a}$ (or on any other coinvariant subspace $\Kth$ of the Hardy space $H^2(\D)$) satisfying \eqref{eq16} is a truncated Toeplitz operator on $K_{\theta_{a}}$. By Theorem~2.1 in \cite{BBK}, for every positive bounded truncated Toeplitz operator $\tilde T$ on $K_{\theta_{a}}$ there exists a finite positive measure $\tilde\mu$ on $\T$ such that $\tilde\mu(\{1\}) = 0$ and
$$
(\tilde T h, h)_{L^2(\T)} = \int_{\T}|h|^2\,d\tilde\mu 
$$    
for all continuous functions $h$ in $K_{\theta_{a}}$. Changing variables in the last integral, we find a positive measure $\mu$ on $\R$ such that
$$
\int_{\T}|h|^2\,d\tilde\mu  = \int_{\R}|f|^2\,d\mu, \qquad f =  Vh.
$$
It follows that $(T f, f) = (\tilde T h, h)_{L^2(\T)} = (A_{\mu, a} f, f)_{L^2(\R)}$ for a dense set of functions $f$ in $\pw_{[0,a]}$. Since $T$ is continuous, we have $T = A_{\mu,a}$. The second part of the Lemma is a direct consequence of relation $T_{\mu, a} = V_{a}^{-1} A_{\mu, 2a} V_{a}$. \qed   

\medskip

\subsection{Preliminaries on canonical Hamiltonian systems.} Let $\H$ be a Hamiltonian on $[0,a]$ with $\trace\H \in L^1[0,a]$. Assume that there is no interval $(r_1, r_2) \subset [0,a]$ such that $\H(t)$ is a constant matrix of rank one for all points~$t \in (r_1, r_2)$. For $r\in[0,a]$ we will denote by $\B(\H,r)$ the de Branges space generated by $\H$ on $[0,r]$, that is, 
$$
\B(\H,r) = \W_{\H,r}L^2(\H,r) = \Bigl\{\mbox{entire }f: \;\; f = \W_{\H,r}X, \;\; X \in L^2(\H,r)\Bigr\},
$$ 
where the Weyl-Titchmarsh transform $\W_{\H,r}$ is defined in \eqref{eq3} for $a = r$. The space $\B(\H,r)$ is actually the Hilbert space with respect to the inner product $(f,g)_{\B(\H,r)} = (f,g)_{L^2(\mu)}$, where $\mu$ is any spectral measure for problem \eqref{eq2}. We refer the reader to paper \cite{BR2015} for the summary of results on direct and inverse spectral theory of canonical Hamiltonian systems and de Brange spaces of entire functions. The readers interested in proofs or in a more detailed account may find necessary information in Chapter\,2 of classical book \cite{dbbook} by L. de Brange or its recent exposition \cite{Romanov} by R.\,Romanov. 

\medskip

\begin{Lem}\label{l4}
Let $\mu$ be an even measure on $\R$ of the form $\mu = cm + \nu$, where $c>0$ and $\nu$ is a finite positive  measure on $\R$ with compact support. Then there exists an infinitely smooth diagonal Hamiltonian $\H$ on $[0,+\infty)$ such that  $\det\H(r) = 1$ for all $r \ge 0$, and $\mu$ is the spectral measure for~$\H$.  
\end{Lem}
\beginpf The result is a kind of folklore. Since the Fourier transform of $\tfrac{1}{c}\nu$ is a smooth (in fact, analytic) function, one can use the classical Gelfand-Levitan approach to find a smooth diagonal potential $Q$ on $[0,a]$ such that $m + \tfrac{1}{c}\nu$ is the spectral measure for the Dirac system $JY' + QY = zY$ corresponding to the boundary condition $Y(0) = \oz$. Then rewrite system $JY' + QY = zY$ as a canonical Hamiltonian system $JX' = z \tilde\H X$ setting $X = M^{-1}Y$, $\tilde\H = M^\ast M$, where $M$ is the matrix solution of equation $JM' = - QM$, $M(0) = \idm$. Observe that $\det\tilde\H = 1$ almost everywhere on $[0,a]$ and $m + \tfrac{1}{c}\nu$ is the spectral measure for system $JX' = z \tilde\H X$, $X(0) = \oz$. To obtain the Hamiltonian on $[0,a]$ corresponding to the spectral measure $\mu$, put $\H = \left(\begin{smallmatrix}\\c & 0\\ 0 &\frac{1}{c}\end{smallmatrix}\right) \tilde\H$. Another (in a sense, equivalent) way of proving Lemma~\ref{l4} is the application of Theorem~5.1 from \cite{Win15}. \qed

\medskip

Define $\type\B(\H,r) = \sup\{\type(f),\, f \in \B(\H,r)\}$ to be the maximal exponential type of entire functions in de Branges space $\B(\H,r)$. The following remarkable formula of Krein \cite{Kr54} and de Brange (Theorem X in \cite{dbii})
\begin{equation}\label{eq17}
\type\B(\H,r) = \int_{0}^{r}\sqrt{\det\H(t)}\,dt,
\end{equation}
represents the maximal exponential type of functions in $\B(\H,r)$ in terms of the Hamiltonian $\H$. Section 6 in \cite{Romanov} contains an elegant self-contained proof of this result. 

\medskip

\begin{Lem}\label{l5}
Let $\H$ be a Hamiltonian on an interval $[0,a]$ such that its spectral measure $\mu$ satisfies \eqref{eq1}. Assume that $\det\H(r) = 1$ for almost all $r\in [0,a]$. Then for all $r\in[0,a]$ we have $\B(\H,r) = (\pw_r,\mu)$. 
\end{Lem}
\beginpf Let $r\in [0,a)$ and let $\eps>0$ be such that $r\in [\eps, a-\eps)$. Then the Hilbert space $(\pw_{r+\eps},\mu)$ of entire functions satisfies an axiomatic description of de Branges spaces (Theorem 23 in \cite{dbbook}) and the embedding $(\pw_{r+\eps},\mu) \subset L^2(\mu)$ is isometric. Since $\mu$ is a spectral measure for $\H$, the embedding $\B(\H,r) \subset L^2(\mu)$ is isometric as well. Applying de Branges chain theorem (Theorem 35 in \cite{dbbook}), we see that ether $(\pw_{r+\eps},\mu) \subset \B(\H,r)$ or $\B(\H,r) \subset (\pw_{r+\eps},\mu)$. Since $\det\H = 1$ almost everywhere on $[0,a]$, formula \eqref{eq17} implies the second alternative. Analogously, one can show that $(\pw_{r-\eps},\mu) \subset \B(\H,r)$. Since this holds for every small number $\eps$ and $\mu$ is sampling, we have $\B(\H,r) = (\pw_{r},\mu)$. 
Finally, for $r = a$ we have 
$$
\B(\H,a) = \ov{\bigcup_{0<r<a} \B(\H,r)} = (\pw_a,\mu),
$$
where the completion is taken with respect to the norm inherited from $L^2(\mu)$. \qed

\medskip

Let $\Theta_\H$ be the absolutely continuous solution of Cauchy problem \eqref{eq2} on $[0,a]$, and denote $\Theta_{\H}^{+} = \langle\Theta_\H, \oz\rangle$, $\Theta_{\H}^{-} = \langle\Theta_\H, \zo\rangle$. The reproducing kernel $k_{\B(\H, r); \lambda}$ at a point $\lambda \in \C$ of the Hilbert space of entire functions $\B(\H,r)$ has the the form
\begin{equation}\label{eq18}
k_{\B(\H, r); \lambda} = \frac{1}{\pi}
\frac{\Theta_{\H}^{+}(r,z)\Theta_{\H}^{-}(r, \bar \lambda) - \Theta_{\H}^{-}(r,z)\Theta_{\H}^{+}(r,\bar \lambda)}{z - \bar \lambda}, \qquad z \in \C.
\end{equation} 
The Paley-Wiener space $\pw_r$ is the de Branges space $\B(\H_0,r)$ for the Hamiltonian $\H_0 = \idm$. The reproducing kernel of $\pw_r$ at $\lambda \in \C$ will be denoted by $\sinc_{r,\lambda}$:
$$
\sinc_{r,\lambda} = \frac{\sin r(z-\bar\lambda)}{\pi(z-\bar\lambda)}, \qquad z \in \C.
$$ 
Using integration by parts and equation \eqref{eq2}, it is easy to show that for each $\lambda \in \C$ we have
$$
\W_{\H,r}\Theta_{\H}(\cdot, \bar\lambda) = \sqrt{\pi} k_{\B(\H, r);\lambda}, \qquad \W_{\H_0,r}\Theta_{\H_0}(\cdot, \bar\lambda) = \sqrt{\pi}\sinc_{r,\lambda},
$$ 
where $\Theta_{\H}(\cdot, \bar\lambda)$ denotes the mapping $t \mapsto \Theta_{\H}(t, \bar\lambda)$ and $\Theta_{\H_0}(\cdot, \bar\lambda)$ is defined analogously.

\medskip

Next assertion is Lemma~4.2 in  \cite{BR2015}.  

\begin{Lem}\label{l6}
Let $\mu$ be a sampling measure for $\pw_a$ and let $r \in [0,a]$. The reproducing kernel of the space $(\pw_r,\mu)$ at $\lambda \in \C$ equals 
$T_{\mu, r}^{-1} \sinc_{r,\lambda}$. 
\end{Lem}
\beginpf For every function $f$ in $(\pw_r,\mu) \subset \pw_r$ and every $\lambda \in \C$ we have
$$
f(\lambda) = (f, \sinc_{a,\lambda})_{L^2(\R)} = (f, T_{\mu,r}^{-1}\sinc_{r,\lambda})_{L^2(\mu)},
$$
where we used the fact that $c_1 I \le T_{\mu,a} \le c_2 I$ on $\pw_a$ and hence $T_{\mu,r}$ is bounded and invertible on $\pw_r$. \qed

\medskip

\begin{Lem}\label{l7}
Let $\phi$ be a function on $[0,a]$ such that $e^{|\phi|} \in L^1[0,a]$. Assume that a spectral measure $\mu$ of problem \eqref{eq2} for the canonical Hamiltonian system generated by $\H= \left(\begin{smallmatrix}e^{\phi}&0\\0&e^{-\phi}\end{smallmatrix}\right)$ satisfies \eqref{eq1} for some constants $c_1$, $c_2$. Then function $w = e^\phi$  belongs to the Muckenhoupt class $A_2[0,a]$ and $\|w\|_{A_2[0,2]} \le 2^{28}c^{14}$, where $c = c_{1}^{-1} + c_{2}^{2} c_{1}^{-1}$. We also have $\frac{1}{a}\int_{0}^{a}(w + \frac{1}{w})\,dx \le 4c$.
\end{Lem}
\beginpf Let us obtain estimates \eqref{eq6}, \eqref{eq7}  for the function $\phi$ as it was suggested in 
Proposition 3.2 of \cite{BR2015}. Take $r \in [0,a]$.
Set $\H_0 = \idm$ and consider the corresponding Weyl-Titchmarsch transforms 
$$
\W_{\H_0,r}: L^2(\H_0,r) \to \B(\H_0,r), \qquad \W_{\H,r}: L^2(\H,r) \to \B(\H,r).
$$ 
We have $\B(\H_0,r) = \pw_r$ and $\B(\H,r) = (\pw_r,\mu)$, see Lemma \ref{l5}.
Since $\mu$ satisfies~$\eqref{eq1}$, the spaces $\pw_r$, $(\pw_r,\mu)$ coincide as sets and 
$$
c_{2}^{-1}\|f\|_{L^2(\R)}^{2} \le \|T_{\mu,r}^{-1} f\|_{L^2(\mu)}^{2} \le c_{1}^{-1}\|f\|_{L^2(\R)}^{2}
$$ 
for every function $f \in \pw_r$. Hence, the operator $T = \W_{\H,r}^{-1}T_{\mu,r}^{-1}\W_{\H_0,r}$ from $L^2(\H_0,r)$ to $L^2(\H,r)$ is correctly defined, bounded, and invertible. 
Moreover, 
\begin{equation}\label{eq19}
c_{2}^{-1}\|X\|_{L^2(\H_0,r)}^{2} \le \|TX\|_{L^2(\H,r)}^{2} \le c_{1}^{-1}\|X\|_{L^2(\H_0,r)}^{2}
\end{equation} 
for every $X \in L^2(\H_0)$. Next, by Lemma \ref{l6} for each $z \in \C$ we have 
$$
T\Theta_{\H_0}(\cdot, z) = \W_{\H,r}^{-1}\bigl(\sqrt{\pi}T_{\mu,r}^{-1}\sinc_{r,\bar z}\bigr) = \Theta_{\H}(\cdot, z). 
$$
For $z = 0$ and all $t \in [0,r]$ we have $\Theta_{\H}(t, 0) = \Theta_{\H_0}(t, 0) = \oz$, hence 
$$
c_{2}^{-1}\|\oz\|_{L^2(\H_0,r)}^{2} \le \|\oz\|_{L^2(\H,r)}^{2} \le c_{1}^{-1}\|\oz\|_{L^2(\H_0,r)}^{2}.
$$
This relation is inequality \eqref{eq7} for the function $\phi$ and constants $b_1 = c_{2}^{-1}$, $b_2 = c_{1}^{-1}$.

\medskip

Now let $\partial_0^n\Theta_{\H}(\cdot, 0)$ denote the derivative of order $n$ of the mapping $z \mapsto \Theta_{\H}(\cdot, z)$ from $\C$ to $L^2(\H,r)$  at the point $z = 0$. Then $T\partial_0^n\Theta_{\H}(\cdot, 0) = \partial_0^n\Theta_{\H_0}(\cdot, 0)$ for all integers $n\ge1$.
The right inequality in \eqref{eq19} yields 
\begin{equation}\label{eq21}
\|\partial_{0}^{n}\Theta_{\H}(\cdot, 0)\|_{L^2(\H,r)}^{2} \le c_{1}^{-1} \|\partial_{0}^{n}\Theta_{\H_0}(\cdot, 0)\|_{L^2(\H_0,r)}^{2}.
\end{equation}
From equation \eqref{eq2} we obtain
\begin{gather}
\partial^{n}_{0}\Theta_{\H}(t, 0) 
= n!\int_{0}^{t}\!\!\int_{0}^{t_1}\!\!\dots\int_{0}^{t_{n-1}}\!J^*\H(t_1)J^*\H(t_2)\ldots J^*\H(t_n)\oz d t_n \ldots d t_1, \label{eq22} \\
\partial^{n}_{0}\Theta_{\H_0}(t, 0) \label{eq23}
= {J^\ast}^{n}\left(\!\begin{smallmatrix}t^n\\0\end{smallmatrix}\!\right), 
\end{gather}
for all $t \in [0,r]$ and $n \ge 1$. Observe that
$$
J^*\H(t_1)J^*\H(t_2)\ldots J^*\H(t_n)\oz =
\begin{cases}
\left(\begin{smallmatrix}0\\(-1)^{\frac{n+3}{2}}\exp(G_{\phi,n}(t))\end{smallmatrix}\right), &\mbox{$n$ is odd,} \\[0.9em]
\left(\begin{smallmatrix}(-1)^{\frac{n}{2}}\exp(G_{\phi,n}(t))\\0\end{smallmatrix}\right), &\mbox{$n$ is even,}
\end{cases} 
$$
where $t = (t_1, \ldots, t_n)$ is a point in simplex $K_{t,n}$, and $G_{\phi,n}$ is defined on $K_{t,n}$ by formula \eqref{eq5}. 
Substitute this representation of $J^*\H(t_1)J^*\H(t_2)\ldots J^*\H(t_n)\oz$ to \eqref{eq22}. Then  \eqref{eq22}, \eqref{eq23}, and  \eqref{eq21} give us inequality \eqref{eq6} for all $n \ge 1$ and all $r \in [0,a]$. It remains to use Proposition \ref{p1} to see that $w \in A_{2}[0,a]$ and $\|w\|_{A_2[0,a]} \le 2^{28}c^{14}$. The estimate 
$\frac{1}{a}\int_{0}^{a}(w + \frac{1}{w})\,dx \le 4c$ follows from Lemma \ref{l1}. \qed

\bigskip

\subsection{Proof of Theorem \ref{t1}.} 
Let $\mu$ be a measure on $\R$ such that estimate \eqref{eq1} holds for some $a>0$. Consider the truncated Toeplitz operator $T_{\mu, a} = T_{\mu}$ on $\pw_{a}$. We have $c_1 I \le T_{\mu} \le c_2 I$, where $I$ stands for the identity operator on $\pw_a$. The operator $T_\mu - c_1 I$ satisfies assumptions of Lemma \ref{l3}. Hence, there exists a measure $\nu \ge 0$ on $\R$ such that $T_{\nu} = T_\mu - c_1 I$. One can suppose that $\nu$ is even (otherwise consider the measure $\tilde\nu$ such that $\tilde\nu(S) = \frac{1}{2}(\nu(S) + \nu(-S))$, and note that $T_{\nu} = T_{\tilde\nu}$). Define a sequence of measures $\mu_j$ by $\mu_j = c_1 m + \chi_{j}\nu$, where
$m$ is the Lebesgue measure on $\R$, and $\chi_{j}$ denotes the indicator function of the interval $[-j,j]$. For every $j \ge 1$ the measure $\mu_j$ is even and satisfies relation $\eqref{eq1}$ with the same constants $c_1$, $c_2$. Indeed, $T_{\mu_j} = c_1 I + T_{\chi_j \nu}$ and
$$
c_1I \le  c_1I + T_{\chi_j \nu} \le c_1I + T_{\nu} = T_{\mu} \le c_{2}I.
$$
By Lemma \ref{l4} and Lemma \ref{l7},  for every $j$ there exists a smooth function $w_j>0$ on the interval $[0, a]$ such that $\|w_j\|_{A_2[0,a]} \le 2^{28}c^{14}$, $c = c_{1}^{-1} + c_{2}^{2}c_{1}^{-1}$, and $\mu_j$ is the spectral measure for the Hamiltonian $\H_j = \left(\begin{smallmatrix}w_j &0 \\ 0 & \frac{1}{w_{j}}\end{smallmatrix}\!\right)$ on $[0,a]$. We also have $\frac{1}{a}\int_{0}^{a}(w + \frac{1}{w})dx \le 4c$ for all $j \ge 1$. This allows us to use ``a reverse H\"older inequality'' for weights in $A_2[0,a]$. It says that for every $C_1>0$ there exist $p>1$ and $C_2>0$ such that for all $h \in A_2[0,a]$ with $\|h\|_{A_2[0,a]} \le C_1$ we have  
$$
\frac{1}{a}\int_{0}^{a} h(x)^{p}\,dx \le C_2\left(\frac{1}{a}\int_{0}^{a} h(x)\,dx\right)^{p}. 
$$
Explicit relations between $C_1$, $C_2$, and $p$ can be found in \cite{Vas03}. From here we see that sequences $\{w_j\}_{j\ge1}$, $\bigl\{\frac{1}{w_j}\bigr\}_{j\ge1}$ are informly bounded in $L^p[0,a]$ for some $p>1$. Hence we can find subsequences $w_{j_k}$, $w^{-1}_{j_k}$ converging weakly in $L^p[0,a]$ to functions $w$, $v$, correspondingly. To simplify notations, let the sequences $\{w_j\}_{j\ge1}$, $\bigl\{\frac{1}{w_j}\bigr\}_{j\ge1}$ themselves be weakly convergent. Let us show that $v = w^{-1}$ almost everywhere on the interval $[0,a]$. This is not always the case for arbitrary weakly convergent sequences in $L^p[0,a]$. 

\medskip

For $z \in \C$ denote by $\Theta_j(\cdot, z)$ solution of equation \eqref{eq2} for the Hamiltonian $\H_j$. Integrating \eqref{eq2}, we get 
\begin{equation}\label{eq24}
J\Theta_j(r, z) - \oz = z\int_{0}^{r}\H_j(t) \Theta_j(t, z)\,dt. 
\end{equation}  
Then for every $j \ge 1$ and $r, r' \in [0,a]$ we have the estimates
\begin{align*}
\|\Theta_{j}(r,z)\|_{\C^2} &\le \exp\left(|z|\int_{0}^{a}\|\H_j(t)\|\,dt\right),\\
\|\Theta_{j}(r,z) - \Theta_{j}(r',z)\|_{\C^2} &\le |z|\cdot |r - r'|^{\frac{p-1}{p}}\left(\int_{0}^{a}\|\H_j(t)\|^p\cdot\|\Theta_{j}(t,z)\|^{p}_{\C^2}\,dt\right)^{\frac{1}{p}},
\end{align*}
showing that functions $\Theta_{j}(\cdot,z)$ are uniformly bounded and equicontinuous on $[0,a]$. Therefore, there is a subsequence of the sequence $\Theta_{j}(\cdot,z)$ converging uniformly on $[0,a]$ to a function $\Theta(\cdot, z)$. As before, we suppose that the sequence $\Theta_{j}(\cdot,z)$ itself is uniformly convergent on $[0,a]$. It is clear that the limit function $\Theta$ satisfies equation~\eqref{eq24} for the Hamiltonian $\H = \left(\begin{smallmatrix}w &0\\0&v\end{smallmatrix}\right)$. Hence, it satisfies equation \eqref{eq2} for~$\H$. Fix a number $r \in (0,a]$. For every $\lambda$ and $z$ in $\C$ we have
\begin{equation}\label{eq25}
k_{\B(\H,r);\lambda}(z) = \lim_{j \to \infty} k_{\B(\H_j,r);\lambda}(z) = \lim_{j \to \infty} (T_{\mu_{j},r}^{-1}\sinc_{r,\lambda})(z) = (T_{\mu,r}^{-1}\sinc_{r,\lambda})(z).  
\end{equation}
Indeed, the first equality above follows from formula \eqref{eq18} and convergence of $\Theta_j$ to $\Theta$ on $[0,a]$ when a spectral parameter ($\bar\lambda$ or $z$) is fixed. Lemma \ref{l5} and Lemma~\ref{l6} give us the second equality. Finally, using the fact that the operators $T_{\mu_{j},r}$ on $\pw_r$ tend to $T_{\mu,r}$ in the strong operator topology, we obtain the last equality in \eqref{eq25}. 
From \eqref{eq25} we see  that Hilbert spaces of entire functions $\B(\H,r)$, $(\pw_{r},\mu)$ have the same reproducing kernels. Hence $\B(\H,r) = (\pw_{r},\mu)$ and formula \eqref{eq17} implies 
$$
r = \int_{0}^{r}\sqrt{\det\H(t)}\,dt, \qquad r \in [0,a].
$$
It follows that $\det\H=1$ almost everywhere on $[0,a]$, that is, $v = w^{-1}$. Next, from the direct spectral theory we know that the family $\{\Theta(\cdot, \lambda)\}_{\lambda \in \C}$ is complete in $L^2(\H,a)$ and $\W_{\H,a} \Theta(\cdot, \lambda) = k_{\B(\H,a);\lambda}$ for every $\lambda \in \C$, where $\W_{\H, a}$ denotes the Weyl-Titchmarsch transform associated to $\H$. Using \eqref{eq25} again, we get
\begin{multline*}
(\Theta(\cdot, \lambda),\Theta(\cdot, z))_{L^2(\H,a)} = \pi k_{\B(\H,a);\lambda}(z) = \pi (T_{\mu,a}^{-1}\sinc_{a,\lambda}, \sinc_{a,z})_{L^2(\R)} \\
= \pi (T_{\mu,a}^{-1}\sinc_{a,\lambda}, T_{\mu,a}^{-1}\sinc_{a,z})_{L^2(\mu)} = (\W_{\H, a}\Theta(\cdot, \lambda),\W_{\H, a}\Theta(\cdot, z))_{L^2(\mu)}.
\end{multline*}
Hence, the operator $\W_{\H, a}$ acts isometrically from $L^2(\H; a)$ to $L^2(\mu)$ and $\mu$ is a spectral measure for 
$\H$. In particular, we can apply Lemma \ref{l7} to $\H$, $\mu$, and conclude that the function $w = e^{\phi}$ is in $A_2[0,a]$ and $\|w\|_{A_2[0,a]} \le 2^{28}c^{14}$. Uniqueness of the Hamiltonian $\H$  follows immediately from formula \eqref{eq25}:
\begin{equation}\label{eq20}
\int_{0}^{r} w(t)\,dt = \int_{0}^{r}\langle\H(t)\oz,\oz\rangle\,dt =  \pi k_{\B(\H,a);0}(0) = \pi \|T_{\mu,r}^{-1}\sinc_{r,0}\|^{2}_{L^2(\mu)},
\end{equation}
where the right hand side is completely determined by $\mu$, while the left hand side determines $\H$.
\qed   

\medskip

\vbox{
Differentiating formula \eqref{eq20}, we obtain the following corollary.
\begin{Cor}\label{c1}
The Hamiltonian $\H = \left(\begin{smallmatrix}w &0\\0&\frac{1}{w}\end{smallmatrix}\right)$ in Theorem \ref{t1} could be recovered from~$\mu$ by means of the following formula: $w(r) = \pi \frac{\partial}{\partial r} \|T_{\mu,r}^{-1}\sinc_{r,0}\|^{2}_{L^2(\mu)}$, $r \in [0,a]$.
\end{Cor}
}

\section{Proof of Theorem \ref{t2} and Theorem \ref{t3}}

\medskip

Let us first show that Theorem \ref{t2} does not follow from a general theory of canonical Hamiltonian systems. Consider the simplest case where the Hamiltonian $\H$ coincides with the identity matrix $\idm$ on $[0,a]$. We claim that there is no unitary operator $U : L^2(\H, a) \to \pw_{[0,2a]}$ such that $UL^2(\H, r) = \pw_{[0,2r]}$ for all $r \in [0,a]$. Indeed, existence of such a unitary operator yields the existence of another unitary operator $\tilde U: L^2[-a,a] \to L^2[0,2a]$ such that $\tilde U L^2[-r,r] = L^2[0,2r]$ for all 
$r \in [0,a]$. For every $r_1 > r_2 \ge 0$ let $\chi_{[r_1, r_2]}$ denote  the indicator function of the interval $[r_1,r_2]$. Put $g = \tilde U \chi_{[0,a]}$ and consider decomposition $g = f_{r} + h_{r}$, where $f_r = \tilde U \chi_{[0,r]}$, $h_r = \tilde U \chi_{[r,a]}$, $r \in [0,a]$. Since $\tilde U L^2[-r,r] = L^2[0,2r]$ by our assumption, the function $f_r$ is supported on $[0,2r]$. Note also that the function $h_r$ is orthogonal to all functions from $L^2[0,2r]$ and hence it is supported on $[2r,2a]$. From here we see that $f_r = \chi_{[0,2r]}g$ for all $r \in [0,a]$. Next, unitarity of the operator $\tilde U$ implies that
$$
\int_{0}^{2r}|g(t)|^{2}\,dt = \int_{0}^{2r}|f_r(t)|^{2}\,dt = \int_{-a}^{a}|\chi_{[0,r]}(t)|^2\,dt = r, \qquad r \in [0,a].
$$   
It follows that $|g(t)|^2 = 1/2$ for almost all $t\in [0,2a]$. In particular, the linear span of functions $f_r \in \tilde U L^2[0,a]$, $r \in [0,a]$, is dense in $L^2[0,2a]$. This contradicts to the fact that $\tilde U$ is a unitary operator from $L^2[-a,a]$ to $L^2[0,2a]$. Thus, the Weyl-Titchmarsh transform $\W_{\H,a}$ from formula \eqref{eq3} can not be used to construct the operator $\F_{\mu}$  from Theorem \ref{t2} by means of superpositon with some simple unitary operators like shifts, reflections, etc. 

\medskip

The main point that helps in proof of Theorem \ref{t2} is the fact that Hamiltonian~$\H$ generated by an even sampling measure for the Paley-Wiener space $\pw_a$ must have rank two almost everywhere on its domain of definition. It is an open question if this is true for general (not necessarily even) sampling measures for $\pw_a$. See also Proposition \ref{p2} in Section \ref{s5} for more details. 

\medskip

\noindent{\bf Proof of Theorem \ref{t2}.} Fix an even sampling measure $\mu$ and construct the Hamiltonians $\H_j$, $\H$, on $[0,a]$ as in the proof of Theorem \ref{t1}. Put $\phi_j = \log w_j$ and $\phi = \log w$, where $w_j$, $w$ are the functions generating $\H_j$, $\H$. Recall that $w_{j}$ tend to $w$ weakly in $L^p[0,a]$ for some $p>1$ and the same is true for $w_{j}^{-1}$ and $w^{-1}$. Let $\Theta_j$, $\Theta$ be the solutions of system \eqref{eq2} generated by Hamiltonians $\H_j$, $\H$, correspondingly. As we have seen, the functions $\Theta_j(\cdot ,z) = \Bigl(\begin{smallmatrix}\Theta_j^+ \\ \Theta_j^- \end{smallmatrix}\Bigr)$ converge uniformly to $\Theta(\cdot, z) = \Bigl(\begin{smallmatrix}\Theta^+ \\ \Theta^- \end{smallmatrix}\Bigr)$ on the interval $[0,a]$ when $z \in \C$ is fixed.   
For $r \in [0, a]$, define entire functions $P_{2r,j}$ and $P_{2r,j}^{*}$ by
\begin{align*}
&P_{2r,j}\colon z \mapsto e^{irz}\left(e^{\frac{\phi_j(r)}{2}}\Theta_{j}^{+}(r,z) - i e^{-\frac{\phi_j(r)}{2}} \Theta_{j}^{-}(r,z)\right), \\ 
&P_{2r,j}^{*}\colon z \mapsto e^{irz}\left(e^{\frac{\phi_j(r)}{2}}\Theta_{j}^{+}(r,z) + i e^{-\frac{\phi_j(r)}{2}} \Theta_{j}^{-}(r,z)\right),
\end{align*}
and let $P_{2r}$, $P_{2r}^{*}$ be defined similarly with $\phi_j$ replaced by $\phi$. These functions satisfy the Krein system of differential equations:
\begin{equation}\label{eq26}
\begin{cases}
P'_{r,j}(z) = iz P_{r,j}(z) + \frac{\phi'_{j}(r/2)}{4}P_{r,j}^{*}(z), &P_{0,j}(z) = e^{\frac{\phi_j(0)}{2}},\\
P^{*}_{r,j}\!\!\!'\,\,(z) = \frac{\phi'_{j}(r/2)}{4} P_{r,j}(z), &P_{0,j}^{*}(z) = e^{\frac{\phi_j(0)}{2}},
\end{cases}
\end{equation}
where $\phi'_{j}(r/2)$ is the value of smooth function $\phi'_j$ at $r/2$. From system \eqref{eq26} we obtain by integration by parts (see Lemma 9.1 in \cite{Den06}) the Christoffel-Darboux formula:
$$
\int_{0}^{r}P_{t,j}(z)\ov{P_{t,j}(\lambda)}\,dt = i \frac{P_{r,j}^{*}(z)\ov{P_{r,j}^{*}(\lambda)} - P_{r,j}(z)\ov{P_{r,j}(\lambda)}}{z - \bar \lambda}.
$$
The right hand side could be rewritten in the form 
$$
\ldots=2 e^{i\frac{r}{2}(z - \bar\lambda)} \cdot \frac{\Theta^{+}_{j}(\frac{r}{2},z)\ov{\Theta^{-}_{j}(\frac{r}{2},\lambda)} - \Theta^{-}_{j}(\frac{r}{2},z)\ov{\Theta^{+}_{j}(\frac{r}{2},\lambda)}}{z - \bar \lambda}, 
$$
which tends to $2\pi k_{r,\lambda}(z)$, the scalar multiple of the reproducing kernel $k_{r,\lambda}$ at $\lambda$ of the Hilbert space 
$e^{i\frac{r}{2}z}\B(\H,\frac{r}{2}) = (\pw_{[0,r]}, \mu)$, see formula \eqref{eq18}. On the other hand, for every pair~$z, \lambda \in \C$ we have 
\begin{align*}
P_{t,j}(z)\ov{P_{t,j}(\lambda)} 
= e^{i\frac{t}{2}(z- \bar\lambda)}\Bigl(&e^{\phi_j(\frac{t}{2})}\Theta^{+}_{j}(\tfrac{t}{2},z)\ov{\Theta^{+}_{j}(\tfrac{t}{2},\lambda)}
+ e^{-\phi_j(\frac{t}{2})}\Theta^{-}_{j}(\tfrac{t}{2},z)\ov{\Theta^{-}_{j}(\tfrac{t}{2},\lambda)} + \\
&+  i \Theta^{+}_{j}(\tfrac{t}{2},z)\ov{\Theta^{-}_{j}(\tfrac{t}{2},\lambda)} - i \Theta^{-}_{j}(\tfrac{t}{2},z)\ov{\Theta^{+}_{j}(\tfrac{t}{2},\lambda)} \Bigr).
\end{align*}
Since functions $e^{\phi_j}$, $e^{-\phi_j}$ converge weakly in $L^p[0,a]$ to functions $e^{\phi}$, $e^{-\phi}$, correspondingly, we see that
\begin{equation}\label{eq27}
\int_{0}^{r}P_{t}(z)\ov{P_{t}(\lambda)}\,dt = \lim_{j \to \infty}\int_{0}^{r} P_{t,j}(z)\ov{P_{t,j}(\lambda)}\,dt = 2\pi k_{r,\lambda}(z)
\end{equation}
for every $r \in [0,2a]$. Let $\chi_r$ be the indicator function of the interval $[0,r]$. Denote by $\Lt$ the set of all finite linear combinations of functions $t \mapsto \chi_r(t)\ov{P_{t}(z)}$ on $[0, 2a]$, where $z \in \C$ and $r \in [0,2a]$. The linear manifold $\Lt$ is dense in $L^2[0,2a]$. Indeed, for every function $g \in L^2[0,2a]$ orthogonal to $\Lt$ we have
$$
0 = \int_{0}^{2a} g(t)\chi_r(t)P_t(0)\,dt = \int_{0}^{r}g(t)e^{\frac{\phi(t/2)}{2}}\,dt, \quad r \in [0,2a], 
$$
yielding $g = 0$ in $L^2[0,2a]$. Formula \eqref{eq27} also shows that a nontrivial finite linear combination of functions $\chi_r(t)\ov{P_{t}(z)}$ cannot vanish almost everywhere on $[0,2a]$. Consider the operator $\F_\mu: L^2[0, 2a] \to (\pw_{[0,2a]},\mu)$ densely defined on $\Lt$ by
$$
\F_\mu: f \mapsto \frac{1}{\sqrt{2\pi}}\int_{0}^{2a}f(t) P_{t}(z)\,dt, \qquad z \in \C.
$$
The operator $\F_\mu$ takes the function $t \mapsto \chi_r(t)\ov{P_{t}(\lambda)}$ on $[0, 2a]$ into $\sqrt{2\pi} k_{r,\lambda}$, see formula \eqref{eq27}. Moreover, for every $r_1, r_2 \in [0,2a]$ we have 
\begin{align*}
\Bigl(\F_\mu \chi_{r_1}\ov{P_{t}(\lambda)}, \F_\mu \chi_{r_2}\ov{P_{t}(z)}\Bigr)_{L^2(\mu)}
&= 2\pi (k_{r_1, \lambda}, k_{r_2, z})_{L^2(\mu)} = 2\pi k_{r, \lambda}(z), \\
&= \Bigl(\chi_{r_1} \ov{P_{t}(\lambda)}, \chi_{r_2}\ov{P_{t}(z)}\Bigr)_{L^2[0,2a]},
\end{align*}
where $r = \min(r_1, r_2)$. This shows that $\F_\mu$ is an isometry on $\Lt$. Since the linear span of the set $\{k_{2a,\lambda}, \; \lambda \in \C\}$ is complete in $(\pw_{[0,2a]}, \mu)$, the operator $\F_\mu$ is unitary. It is also clear from the definition that $\F_\mu$ maps $L^2[0,r]$ onto $(\pw_{[0,r]},\mu)$ for every $r \in [0,2a]$.\qed

\medskip

\noindent{\bf Proof of Theorem \ref{t3}.} At first, consider a positive bounded invertible operator $W_\psi$ with real symbol $\psi \in \Sch'$ on a finite interval $[0,a]$. Let $\F$ denote the unitary Fourier transform on $L^2(\R)$. Take a smooth function $h$ with support in $(0,a)$ and put $\hat f = \F f$. Consider the operator $\hat{W}_{\psi} = \F W_{\psi} \F^{-1}$ on $\pw_{[0,a]}$. We have
$(\hat{W}_{\psi} \hat h, \hat h)_{L^2(\R)} =  \bigl\langle\hat{\psi}, |\hat{h}|^2\bigr\rangle_{\Sch'}$,
where $\hat\psi$ is the Fourier transform of the tempered distribution $\psi$. It follows that 
$$(\hat{W}_{\psi} f, f)_{L^2(\R)} = (\hat{W}_{\psi} \tfrac{z-i}{z+i}f, \tfrac{z-i}{z+i}f)_{L^2(\R)}$$ 
on a dense subset of the set 
$Z_{-i} = \{f \in \pw_{[0,a]}:\; f(-i) = 0\}$. Since $\hat{W}_{\psi}$ is bounded on $\pw_{[0,a]}$, 
we have the last identity for all $f \in Z_{-i}$.    
Hence, the operator $\hat{W}_{\psi}$ satisfies assumptions of Lemma \ref{l3} and we can find a positive Borel measure~$\mu$ on~$\R$ such that 
$$
(\hat{W}_{\psi} f, g)_{L^2(\R)} = \int_{\R}f \bar g\,d\mu
$$ 
for all $f,g \in \pw_{[0,a]}$. As in the proof of Theorem \ref{t1}, we can assume that the measure $\mu$ is even. Indeed, since $\psi$ is real, we have $(\hat{W}_{\psi} f, f) = (\hat{W}_{\psi} f^*, f^*)$ for arbitrary $f \in \pw_{[0,a]}$  and its reflection $f^*: x \mapsto f(-x)$. By the assumption, the operator $\hat{W}_{\psi}$ is positive, bounded and invertible on $\pw_{[0,a]}$. Hence the measure $\mu$ satisfies \eqref{eq1} for some $c_1, c_2$ and $a/2$ in place of $a$. By Theorem \ref{t2}, there is a unitary operator $\F_\mu: L^2[0,a] \to (\pw_{[0,a]}, \mu)$ such that $\F_\mu: L^2[0,r] = (\pw_{[0,r]}, \mu)$ for every $r \in [0,a]$. Identifying Hilbert spaces $(\pw_{[0,a]}, \mu)$ and $\pw_{[0,a]}$ as sets, we can define the operator $A = \F_{\mu}^{-1}\F$  on $L^2[0,a]$. By construction, the operator $A$ is bounded and invertible  and $A L^2[0,r] = L^2[0,r]$ for every $r \in [0,a]$. We also have
\begin{equation}\label{eq33}
(W_{\psi}h,h)_{L^2[0,a]} = \int_{\R} \bigl| \hat h\bigr|^{2}\,d\mu = (\F h, \F h)_{L^2(\mu)} 
= (\F_{\mu}^{-1} \F h, \F_{\mu}^{-1} \F h)_{L^2[0,a]} 
\end{equation}
for all smooth functions $h$ with support in $(0,a)$. It follows that the operator $W_{\psi}$ admits the triangular factorization $W_{\psi} = A^*A$. 

\medskip

It remains to consider the case where $W_\psi$ is a positive bounded invertible Wiener-Hopf operator on $L^2[0,\infty)$ with real symbol $\psi \in \Sch'$. It is known (see Section 4.2.7 in \cite{Nik02t1}) that in this case the Fourier transform of the distribution $\psi$ coincides with a function $\sigma$ on $\R$ such that $c_1 \le \sigma(x) \le c_2$ for some positive constants $c_1$, $c_2$ and almost all $x \in \R$. In particular, the measure $\mu = \sigma\,dm$ is sampling for all Paley-Wiener spaces $\pw_{[0,r]}$, $r > 0$. Since $\psi$ is real, the function $\sigma$ is even. For every $r>1$ we can use Theorem \ref{t1} and find a Hamiltonian $\H_r$ on $[0,r]$ such that $\det\H_r(t) = 1$ for almost all $t \in [0,r]$ and $\mu$ is the spectral measure for $\H_r$. Since the Hamiltonian $\H$ in Theorem \ref{t1} is defined uniquely, we have $\H_{r}(t) = \H_{r'}(t)$ for almost all $t \in [0, \min(r, r')]$. 
This shows that there is the Hamiltonian $\H$ on $[0,\infty)$ such that $\det \H = 1$ almost everywhere and $\mu$ is the spectral measure for $\H$. In particular, we can define a family of entire functions 
$\{P_{t}\}_{t \ge 0}$ such that the mapping
\begin{equation}\label{eq32}
\F_\mu: f \mapsto \frac{1}{\sqrt{2\pi}}\int_{0}^{r}f(t)P_t(z)\,dt
\end{equation}
sends unitarily the space $L^2[0,r]$ onto the space $(\pw_{[0,r]}, \mu)$ for every $r>0$, see the proof of Theorem \ref{t2}. Let $H^2_\mu(\C_+)$ be the weighted Hardy space with the inner product $(f,g)_{H^2_\mu(\C_+)} = (f,g)_{L^2(\mu)}$. Since $c_1 \le \sigma \le c_2$ on $\R$, the space $H^2_\mu(\C_+)$ coincides as a set with the standard Hardy space $H^2(\C_+) = \F L^2[0,\infty)$. Define the unitary operator $\F_\mu$ from $L^2[0,\infty)$ to $H^2_\mu(\C_+)$ by formula \eqref{eq32} with $r = \infty$ on  the dense set of compactly supported bounded functions in $L^2[0,\infty)$. Then the operator $A = \F_{\mu}^{-1}\F$ on $L^2[0,\infty)$ is bounded and invertible. Moreover, $A L^2[0,r] = L^2[0,r]$ for every $r \ge 0$, and $W_\psi = A^*A$, see formula \eqref{eq33}. \qed

\medskip

\noindent {\bf Remark.} It can be shown that positive bounded invertible Wiener-Hopf operators $W_\psi$ on $L^2[0,a)$ with real symbols $\psi \in \Sch'$ admit triangular factorisation in the reverse order, $W_\psi = AA^*$. In the case $a = \infty$ the classical Wiener-Hopf factorization works: one can take $A = \F^{-1}T_{\ov{\phi}_\sigma}\F$, where $T_{\phi_\sigma}$ is the Toeplitz operator on $H^2(\C_+)$ with analytic symbol $\phi_\sigma$ such that $|\phi_\sigma|^2 = \sigma = \F\psi$. If $a>0$ is finite, then we can use Theorem \ref{t3} to find left triangular factorization $W_\psi = \tilde A^* \tilde A$ and then put $A = C_{a}\tilde A C_a$, where $C_a: f \mapsto \ov{f(a-x)}$ is the conjugate-linear isometry on $L^2[0,a]$. Since $C_a W_{\psi} C_a = W_{\psi}$ for the self-adjoint Wiener-Hopf operator $W_\psi$ on $L^2[0,a]$, and $C_a^2 = I$, we have $W_\psi = AA^*$. It is also clear that the operator $A$ is upper-triangular.  

\section{Appendix. Two results by L.~A.~Sakhnovich}\label{s5}
In paper \cite{Sakhnovich12} L.~A.~Sakhnovich proved (see Theorem 4.1 and Remark 4.1 in \cite{Sakhnovich12}) that positive bounded invertible  Wiener-Hopf operator 
\begin{equation}\label{eq30}
T: f \mapsto f -\mu\int_{0}^{\infty}f(t)\frac{\sin\pi(t-x)}{\pi(t -x)}\,dt, \qquad f \in L^2[0,\infty), \quad 0<\mu<1,
\end{equation}
densely defined on $L^2[0,\infty)$ does not admit triangular factorization $T = A^*A$, where a bounded invertible operator $A$ on $L^2[0,\infty)$ is such that $AL^2[0,r] = L^2[0,r]$ for every $r \ge 0$. Clearly, this assertion contradicts Theorem \ref{t3}. Let us point out an error in its proof. 

\medskip

The argument in \cite{Sakhnovich12} crucially uses the following claim. Let $\chi_{[-\pi, \pi]}$ be the indicator function of the interval $[-\pi, \pi]$. Formulas $(4.1)-(4.4)$ in \cite{Sakhnovich12} for $n = 0$ and $a_0 = \pi$ determine the function $\sigma': x \mapsto \frac{1}{2\pi}(1 - \mu\cdot \chi_{[-\pi, \pi]}(x))$ on $\R$. The function 
$$
\Pi(z) = \frac{1}{\sqrt{2\pi}}\exp\left(\frac{1}{2i\pi}\int_{-\infty}^{\infty}\frac{1 + tz}{(z-t)(1+t^2)}\log\sigma'(t)\,dt\right)
$$ 
from formula $(4.10)$ of \cite{Sakhnovich12} (see also formula $(4.12)$ therein) is claimed to satisfy the following relation (formula $(4.18)$ in \cite{Sakhnovich12}): 
$$
\lim_{y \to +0}\Pi(iy) = \sqrt{1 - \mu}.
$$
However, this fact is false. Indeed, we have 
$$\frac{1}{\pi i}\frac{1 + tz}{(z-t)(1+t^2)} = -\frac{1}{\pi i}\left(\frac{1}{t - z} - \frac{t}{1+t^2}\right)$$
and hence $\sqrt{2\pi}\Pi(z)$ is the outer function in $\C_+$ whose absolute value on $\R$ coincides with $(\sigma')^{-1/2}$ almost everywhere on $\R$. Since $(\sigma')^{-1/2}$ is regular (in fact, constant) near the origin, we have 
$$\lim_{y \to +0}\Pi(iy) = \frac{1}{\sqrt{2\pi}}(\sigma')^{-1/2}(0) = \frac{1}{\sqrt{1 - \mu}}.$$
We also would like to note that the last relation agrees well with the first identity in formula $(4.19)$ from \cite{Sakhnovich12}.

\medskip

The second part of this section concerns factorization problem for truncated Toeplitz operators generated by general sampling measures for the space $\pw_a$ not necessarily symmetric with respect to the origin. The result is equivalent to Theorem 4.2 in \cite{Sakhnovich07}. The proof below seems to be a bit more straightforward than the original one, possibly, because we consider the one-dimensional situation.  
\begin{Prop}\label{p2} Let $\H$ be a Hamiltonian on $[0,\ell]$ such that $\int_{0}^{\ell}\trace\H(r)<\infty$, and let $\mu$ be a spectral measure for problem \eqref{eq2}. Set $a = \int_{0}^{\ell}\sqrt{\det\H(r)}\,dr$. Assume that $\mu$ satisfies \eqref{eq1}. The following assertions are equivalent:
\begin{itemize}
	\item[$(a)$] $\det\H >0$ almost everywhere on $[0, \ell]$;
	\item[$(b)$] there exists a unitary operator $V_\mu: \pw_{a} \to (\pw_a,\mu)$ such that for every $r \in [0, a]$ we have $V_\mu\pw_{r} = (\pw_a,\mu)$.
  \item[$(c)$] there exists a bounded invertible operator $A$ on $\pw_a$ such that $T_{\mu,a} = A^*A$ and for every $r \in [0,a]$ we have $A\pw_r = \pw_r$.
\end{itemize}
\end{Prop}
Given a Hamiltonian $\H$ on $[0, \ell]$ such that $a = \int_{0}^{\ell}\sqrt{\det\H(t)}\,dt > 0$, we define continuous from the left function $\xi_\H$ from $[0, a]$ to $[0, \ell]$ by 
$$
r = \int_{0}^{\xi_\H(r)}\sqrt{\det\H(t)}\,dt, \qquad r \in [0,a].
$$
This function is continuous if and only if there are no interval $(r_1,r_2) \subset [0, \ell]$ such that $\det\H(t) = 0$ for almost all $t \in (r_1,r_2)$. The function $\xi_\H$ is absolutely continuous if and only if $\det\H(t) > 0$ for almost all $t \in [0,a]$, see Exercise 13 in Chapter IX of \cite{Natanson}.  

\medskip

\noindent{\bf Proof of Proposition \ref{p2}.} $(a) \Rightarrow (b)$. Since $\det\H >0$ almost everywhere on the interval $[0, \ell]$, the function $\xi = \xi_\H$ is absolutely continuous and 
$$
\xi'(r) = \frac{1}{\sqrt{\det{\H(\xi(r))}}}
$$ 
for almost all $r\in [0,a]$.  Consider the Hamiltonian $\Ht : r \mapsto \xi'(r)\H(\xi(r))$ on the interval $[0,a]$. We have $\det \Ht = 1$  and $\Theta_{\Ht}(r,z) = \Theta_{\H}(\xi(r),z)$ on $[0, a]$. Changing variable in \eqref{eq3}, we see that $\B(\tilde\H, r) = \B(\H, \xi(r))$ for every $r \in [0,a]$, hence $\mu$ is the spectral measure for~$\Ht$. Consider the Weyl-Titchmarsh transforms generated by Hamiltonians $\Ht$ and $\H_0 = \idm$, correspondingly,
$$
\W_{\Ht,a}: L^2(\Ht,a) \to \B(\Ht,a), \qquad \W_{\H_0,a}: L^2(\H_0,a) \to \pw_a.
$$ 
Define the operator $V_\mu: \pw_{a} \to \B(\Ht,a)$ by $V_\mu = \W_{\tilde\H,a} \mathbb{M}_{\tilde\H^{-1/2}}\W_{\H_0,a}^{-1}$,
where $\mathbb{M}_{\Ht^{-1/2}}: L^2(\H_0, a) \to L^2(\tilde\H, a)$ is the multiplication operator by $\tilde\H^{-1/2}$, that is,  $\mathbb{M}_{\Ht^{-1/2}}:  X \mapsto \Ht^{-1/2}X$. Since $\mathbb{M}_{\Ht^{-1/2}}$ is unitary, the operator $V_\mu$ is unitary as well. It is also clear that $V_\mu\pw_{r} = \B(\Ht,r)$ for every $r \in [0, a]$. Using Lemma \ref{l5}, we see that $\B(\Ht,r) = (\pw_r,\mu)$, as required.

\medskip

\noindent $(b) \Rightarrow (a)$. We will show that the function $\xi = \xi_\H$ is absolutely continuous. 
Let $\chi_{r}$ be the indicator function of an interval $[0, r]$. For every $r \in [0,a]$ consider the functions
$X_{\xi(r)} = \chi_{\xi(r)}\oz$, $Y_{\xi(r)} = \chi_{\xi(r)}\zo$ in $L^2(\H, \xi(r))$. A straightforward modification of Lemma~\ref{l5} gives $\B(\H,\xi(r)) = (\pw_r,\mu)$ for all $r \in [0,a]$. Put 
$$
X_{r}^{0} = \W_{\H_0,a}^{-1}V_{\mu}^{-1} \W_{\H,a}X_{\xi(r)}, \quad Y_{r}^{0} = \W_{\H_0,a}^{-1}V_{\mu}^{-1} \W_{\H,a} Y_{\xi(r)}.
$$
Since $V_\mu$ is isometric and $V_\mu\pw_r = (\pw_r,\mu)$, we have $\mathcal{P}_{\mu,r} V_\mu = V_\mu \mathcal{P}_r$, where $\mathcal{P}_r$, $\mathcal{P}_{\mu,r}$ are the orthogonal projections on $\pw_a$, $(\pw_a,\mu)$, with ranges $\pw_r$, $(\pw_r,\mu)$, respectively.  It follows that $X_{r}^{0} = \chi_r X_{a}^{0}$ and 
$Y_{r}^{0} = \chi_r Y_{a}^{0}$. Using the fact that the operators $\W_{\H_0,a}$, $\W_{\H,a}$ are unitary, we obtain
\begin{align}
\int_{0}^{\xi(r)} \trace\H(t)\,dt  \notag
&= \int_{0}^{\xi(r)}\Bigl(\bigl\langle \H(t) \oz, \oz \bigr\rangle + \bigl\langle \H(t) \zo, \zo \bigr\rangle\Bigr)\,dt,\\ \notag
&=\|X_{\xi(r)}\|^{2}_{L^2(\H, \ell)} + \|Y_{\xi(r)}\|^{2}_{L^2(\H, \ell)}, \\ \notag
&= \|X_{r}^{0}\|^{2}_{L^2(\H_0, a)} + \|Y_{r}^{0}\|^{2}_{L^2(\H_0, a)},\\ \notag
&=\|\chi_r X_{a}^{0}\|^{2}_{L^2(\H_0, a)} + \|\chi_r Y_{a}^{0}\|^{2}_{L^2(\H_0, a)},\\ 
&=\int_{0}^{r} \Bigl(\|X_{a}^{0}(t)\|_{\C^2}^{2}+\|Y_{a}^{0}(t)\|_{\C^2}^{2}\Bigr) \,dt. \label{eq28}
\end{align}
The above equalities hold for all $r \in [0,a]$. Let us define the function $\kappa$ on $[0, \ell]$ by
$$
\kappa(s) = \int_{0}^{s}\trace\H(t)\,dt, \qquad s \in [0, \ell]. 
$$
Then $\kappa$ is an absolutely continuous function with positive derivative almost everywhere on $[0,\ell]$, hence the inverse mapping $\kappa^{-1}$ is also absolutely continuous and has positive derivative. On the other hand, formula \eqref{eq28} shows that $\kappa(\xi)$ is an absolutely continuous function. It follows that the superposition $\xi = \kappa^{-1}(\kappa(\xi))$ is absolutely continuous and hence $\det\H > 0$ almost everywhere on $[0,\ell]$. 

\medskip

\noindent $(b) \Rightarrow (c)$. Since $\mu$ satisfies \eqref{eq1}, the identical embedding $j: \pw_{a} \to (\pw_a,\mu)$ is a bounded and invertible operator. Define $A = V_{\mu}^{-1} j$. Then for all $f,g$ in $\pw_a$ we have
\begin{equation}\label{eq29}
(A^*A f,  g) = (V_{\mu}^{-1} j f, V_{\mu}^{-1} j g)_{L^2(\R)} = (jf, jg)_{L^2(\mu)} = \int_{\R}f\bar{g}\,d\mu=(T_{\mu,a} f, g),
\end{equation}
by the unitarity of the operator $V_\mu$. It follows that $T_{\mu,a}= A^*A$. By construction, the operator $A$ is invertible. We  also have $A \pw_{r} = \pw_{r}$ for all $r \in [0,a]$, hence $A$ is upper-triangular.

\medskip

\noindent $(c) \Rightarrow (b)$. Assume that $T_{\mu,a}$ admits a left triangular factorization $T_{\mu,a} = A^*A$. Define the operator  $V_\mu: \pw_{a} \to (\pw_\mu,a)$ by 
$V_\mu = jA^{-1}$, 
where $j$ is the embedding from $\pw_{a}$ to $(\pw_{a}, \mu)$. Then $V_\mu \pw_{r} = (\pw_{r}, \mu)$ for every $r \in [0, a]$ and 
\begin{align*}
(V_{\mu} f, V_{\mu} g)_{L^2(\mu)} 
&= ((A^{-1})^* j^* j A^{-1} f,g)_{L^2(\R)}\\
&= ((A^*)^{-1} T_{\mu,a} A^{-1} f,g)_{L^2(\R)} = (f,g)_{L^2(\R)},  
\end{align*}
where we used the identity $T_{\mu,a} = j^* j$, see \eqref{eq29}. Since $A$ and $j$ are invertible, $V_\mu$ is a unitary operator.
\qed

\bibliographystyle{plain} 
\bibliography{bibfile}

\def\cprime{$'$} \def\cprime{$'$} \def\cprime{$'$}
\begin{thebibliography}{10}

\bibitem{BBK}
Anton Baranov, Roman Bessonov, and Vladimir Kapustin.
\newblock Symbols of truncated {T}oeplitz operators.
\newblock {\em J. Funct. Anal.}, 261(12):3437--3456, 2011.

\bibitem{BR2015}
R.~V. Bessonov and R.~V. Romanov.
\newblock An inverse problem for weighted {P}aley-{W}iener spaces.
\newblock {\em preprint arXiv:1509.08117}, 2015.

\bibitem{dbii}
Louis de~Branges.
\newblock Some {H}ilbert spaces of entire functions. ii.
\newblock {\em Trans. Amer. Math. Soc.}, 99:118–152, 1961.

\bibitem{dbbook}
Louis de~Branges.
\newblock {\em Hilbert spaces of entire functions}.
\newblock Prentice-Hall, Inc., Englewood Cliffs, N.J., 1968.

\bibitem{Den06}
Sergey~A. Denisov.
\newblock Continuous analogs of polynomials orthogonal on the unit circle and
  {K}re\u\i n systems.
\newblock {\em IMRS Int. Math. Res. Surv.}, pages Art. ID 54517, 148, 2006.

\bibitem{GL55}
I.~M. Gelfand and B.~M. Levitan.
\newblock {\em On the determination of a differential equation from its
  spectral function}.
\newblock American Mathematical Society, 1955.

\bibitem{GK70}
I.~C. Gohberg and M.~G. Krein.
\newblock {\em Theory and Applications of {V}olterra Operators in {H}ilbert
  Space}, volume~24 of {\em Translations of Mathematical Monographs}.
\newblock American Mathematical Society, 1970.

\bibitem{Krein74}
I.~S. Kac and M.~G. Krein.
\newblock On the spectral functions of the string.
\newblock {\em Amer. Math. Soc. Transl}, 103(2):19--102, 1974.

\bibitem{Kr54}
M.~G. Krein.
\newblock On a fundamental approximation problem in the theory of extrapolation
  and filtration of stationary random processes.
\newblock {\em Dokl. Acad. Nauk SSSR}, 94:13--16, 1954.

\bibitem{Kr81}
M.~G. Krein.
\newblock Continuous analogues of propositions on polynomials orthogonal on the
  unit circle.
\newblock {\em Dokl. Akad. Nauk SSSR (N.S.)}, 105:637--640, 1955.

\bibitem{Larson85}
David~R. Larson.
\newblock Nest algebras and similarity transformations.
\newblock {\em Annals of Mathematics}, 121(2):409--427, 1985.

\bibitem{March11}
Vladimir~Aleksandrovich Marchenko.
\newblock {\em Sturm-{L}iouville operators and applications}, volume 373.
\newblock American Mathematical Soc., 2011.

\bibitem{Natanson}
I.~P. Natanson.
\newblock {\em Theory of functions of a real variable}, volume~1.
\newblock Ungar, New York, 1964.
\newblock Translated by {L}. {F}. {B}oron, with the editorial collaboration of
  and with annotations by {E}. {H}ewitt.

\bibitem{Nik02t1}
Nikolai~K. Nikolski.
\newblock {\em Operators, functions, and systems: an easy reading. {V}ol. 1},
  volume~92 of {\em Mathematical Surveys and Monographs}.
\newblock American Mathematical Society, Providence, RI, 2002.
\newblock Hardy, Hankel, and Toeplitz, Translated from the French by Andreas
  Hartmann.

\bibitem{Romanov}
Roman Romanov.
\newblock Canonical systems and de {B}ranges spaces.
\newblock {\em preprint arXiv:1408.6022}, 2014.

\bibitem{Sakhnovich07}
Lev~A. Sakhnovich.
\newblock On triangular factorization of positive operators.
\newblock In {\em Recent Advances in Matrix and Operator Theory}, pages
  289--308. Springer, 2007.

\bibitem{Sakhnovich12}
Lev~A. Sakhnovich.
\newblock Effective construction of a class of positive operators in {H}ilbert
  space, which do not admit triangular factorization.
\newblock In {\em Levy Processes, Integral Equations, Statistical Physics:
  Connections and Interactions}, pages 85--99. Springer, 2012.

\bibitem{Sar07}
Donald Sarason.
\newblock Algebraic properties of truncated {T}oeplitz operators.
\newblock {\em Oper. Matrices}, 1(4):491--526, 2007.

\bibitem{Vas03}
V.~I. Vasyunin.
\newblock The sharp constant in the reverse {H}\"older inequality for the
  {M}uckenhoupt weights.
\newblock {\em Algebra i Analiz}, 15(1):73--117, 2003.

\bibitem{Win15}
Henrik Winkler.
\newblock {\em Operator Theory}, chapter Two-Dimensional Hamiltonian Systems,
  pages 1--22.
\newblock Springer Basel, Basel, 2014.

\end{thebibliography}
\enddocument